\documentclass[oneside]{amsart}
\usepackage{amsmath, amssymb}
\usepackage{amsthm}
\newtheorem{theorem}{Theorem}[section]

\newtheorem{lemma}[theorem]{Lemma}
\newtheorem{claim}{Claim}
\def\nbd{neighborhood }
\def\nbds{neighborhoods }

\def\u{{\mathop{\uparrow} }}
\def\d{{\mathop{\downarrow} }}
\def\cl{\mathrm{cl }}
\def\X{\widetilde{X}}
\def\a{\alpha}
\def\<{{< \omega}}
\def\V{\bigvee }
\def\W{\bigwedge }

\title{On completeness of H-closed pospaces}
\author{Tomoo Yokoyama}
\date{\today}
\address{Tokyo University of Marine Science and Technology\\
Faculty of Marine Technology\\
1-2-6 Etchujima, Koto-ku, Tokyo, 108-8477 (Japan)}


\begin{document}

\maketitle
 
\begin{abstract}
We generalized the characterization of H-closedness for linearly ordered pospaces as follows: 
A  pospace $X$ without an infinite antichain   
is an H-closed pospace if and only if 
$X$ is a directed complete and down-complete poset with  
$ \V L \in \cl \d L$  
and $ \W L \in \cl \d L$ for any nonempty chain $L \subseteq  X$. 
\end{abstract}

\section{Introduction and notation}

In \cite{GPR}, 
they gave the following characterization of H-closedness for a
linearly ordered pospace to be H-closed: 
A linearly ordered pospace $X$ is H-closed if and only if 
$X$ is a complete lattice with 
$ \V L \in \cl \d L$ and $ \W L \in \cl \d L$ for any nonempty chain $L \subseteq  X$. 

In the same paper, they also gave a non-H-closed pospace such that 
it is a directed complete and down-complete poset with an infinite antichain and with 
$ \V L \in \cl \d L$ and $ \W L \in \cl \d L$ for any nonempty chain $L$ of it.  
(In fact, they showed that the extension $X$ of an infinite antichain by adding the minimal element and 
equipped with the discrete topology on $X$ is a non-H-closed pospace.) 

In this paper, all topological spaces will be assumed Hausdorff. 
We extend the above characterization of H-closedness for linearly ordered pospaces to one for pospaces without an infinite antichain. 

For a set $X$, denote by $X^\<$ the set of finite subsets of $X$.   
If $A$ is a subset of a topological space $X$, then we denote
the closure of the set $A$ in X by $\cl_X A$ or $\cl A$. By a partial order on a set $X$ we mean a reflexive,
transitive and anti-symmetric binary relation $\leq$ on $X$. 
A set endowed with a partial order is called a partially
ordered set (or poset). 

Recall that a poset with a topology defined on it is called a topological partially ordered space (or pospace) 
if the partial order is a closed subset of $X \times X$. 
A partial order $\leq$ is said to be continuous or closed if $x \nleq y$ in $X$ implies that 
there are open \nbds $U$ and $V$ of $x$ and $y$ respectively such that
$\u U \cap V = \emptyset$ (i.e. $U \cap \d V = \emptyset$). 
A partial order $\leq$ on a topological space $X$ is continuous if and only if 
$(X, \leq )$ is a pospace \cite{W}. 
In any pospace, $\d x$ and $\u x$ are both closed for any element $x$ of it. 

A Hausdorff pospace $X$ is said to be an H-closed pospace if $X$ is a closed
subspace of every Hausdorff pospace in which it is contained. 
Obviously that the notion of H-closedness is a generalization of compactness.

For an element $x $ of a poset $X$, 
$\mathop{\uparrow}  x := \{ y \in X \mid x \leq y \} $ 
(resp. $\mathop{\downarrow}  x := \{ y \in X \mid y \leq x \} $) is  
called the upset (resp. the downset) of $x$.  
For a subset $Y \subseteq X$, 
$\mathop{\uparrow}  Y:= \bigcup_{y \in Y}\mathop{\uparrow}  y $ 
(resp. $\mathop{\downarrow}  Y := \bigcup_{y \in Y}\mathop{\downarrow}  y $) is  
called the upset (resp. the downset) of $Y$.  
For a subset $S$ of a poset, denote by $S^\u$ (resp. $S^\d$) the set of upper (resp. lower) bounds of $S$. 
For elements $x, y$ of a poset, $x \| y$ means that $x$ and $y$ are incomparable. 
For an element $x$ of a poset $X$, denote by $I_x$ the set of incomparable elements for $x$ 
(i.e. $I_x = X -(\u x \cup \d x)$). 
For a subset $A$ of a poset, 
$A$ is said to be a chain if $A$ is  linearly ordered, 
and is said to be an antichain if any distinct elements are incomparable.  
A maximal chain is a chain which is properly contained in no other chain.
The Axiom of Choice implies the existence of maximal chains in any poset.
A subset $D$ of a poset $X$ is (up-)directed (resp. down-directed) if 
every finite subset of $D$ has an upper (resp. lower) bound in $D$. 
A poset $X$ is said to be down-complete (resp. up-complete) 
if each down-directed (resp. up-directed) set $S$ of $X$ has $\bigwedge S$ (resp. $\bigvee S$).   
An up-complete poset is also called a directed complete poset or a dcpo. 
It is well-known that 
a poset $X$ is  directed complete if and only if 
each chain $L$ of $X$ has $\bigvee L$.   

Now we state the main result.

\begin{theorem}\label{thm:001}
Let $X$ be a pospace without an infinite antichain. 
Then $X$ is an H-closed pospace if and only if 
$X$ is directed complete and down-directed such that 
$ \V L \in \cl \d L$  
and $ \W L \in \cl \d L$ for any nonempty chain $L \subseteq  X$. 
\end{theorem}

\section{Proofs}

Note that $\u F \cup \d F = X$ 
for a maximal antichain $F$ of a poset $X$. 
Moreover notice that 
that if $X$ has no infinite antichain, then all subposet and 
all extensions of $X$ by adding finite points have no infinite antichain neither. 

\begin{lemma}\label{lem:001}
Let $X$ be a poset and $x$ a point of $X$. 
Suppose that there is a subset $F$ of $X$ such that 
$F \sqcup \{ x \}$ is a maximal antichain in $X$. 
Let $U := X - ( \u F \cup \d F )$. 
Then $\u U \subseteq \u x \cup \d x$ and 
$\d U \subseteq \u x \cup \d x$. 
\end{lemma}

\begin{proof}
Put 
$X_- := \d (F \sqcup \{ x \})$ and 
$X_+ := \u (F \sqcup \{ x \})$. 
Then $X = X_- \cup X_+$. 
Since $\u U \cap \d F = \emptyset$, we have 
$\u U \cap X_- = \u U \cap \d (F \sqcup \{ x \}) 
= \u U \cap \d  x $ and 
so 
$\u (\u U \cap X_-) \cap X_+ \subseteq \u x $. 
Since $U \cap X_+ \subseteq \u x$, 
we obtain 
$\u U \cap X_+ = 
(\u ( U \cap X_+) \cup \u ( U \cap X_-)) \cap X_+ \subseteq \u x \cup \d x$ and 
so $\u U \subseteq \u x \cup \d x$. 
By symmetry, we obtain  
$\d U \subseteq \u x \cup \d x$. 
\end{proof}

\begin{lemma}\label{lem:03}
Let $X$ be an H-closed pospace without an infinite antichain. 
Then any maximal chain of $X$ is complete. 
\end{lemma}

\begin{proof}
Suppose that
there is a maximal chain $L$ which is non-complete.
Then there is a subset $S$ of $L$ such that 
either $\bigvee_L S$ or $\bigwedge_L S$ does not exist. 
We may assume that $\bigvee_L S$ does not exist. 
If $S^\u = \emptyset$, then 
let $\X := X \sqcup \{ \infty \}$ be the extension with the maximal element $\infty$. 
Define a topology $\tau$ on $\X$ by an open subbase $\tau \cup \{ \X - \d F \mid F \in X^\< \}$. 
Then $X$ is an embedded subspace which is not closed.
We will show that $\X$ is a pospace. 
For any element $x$ of $X$, 
if $x \notin \max X$, then the fact that 
$\bigvee_L S$ does not exist implies that 
there is an element $y > x$ of $X$ and a finite subset $F$ of $X$ such that 
$y \in F$ and $F$ is a maximal antichain. 
Then $U := \X - \u F$ is an open \nbd of $x$ and 
$V := \X - \d F$ is an open \nbd of $\infty$ such that 
$\d U \cap V = \emptyset$. 
Otherwise $x \in \max X$. 
Then there is a finite subset $F$ of $X$ such that  
$F \sqcup  \{ x\}$ is a maximal chain. 
Let $U = X - ( \d _X F \cup \u_X F) \subseteq \d_X x$ be an open \nbd of $x$ 
and 
$V = \X - \d (  F \sqcup \{ x \} )$ an open \nbd of $\infty$. 
Then $\d U \cap V = \emptyset$.  
Thus $\X$ is a pospace and so $X$ is not H-closed, which contradicts. 
Thus $S^\u \neq \emptyset$. 
Then $\W_L(L \setminus \d S)$ does not exist. 
Let $A := \u ( L \setminus \d S)$ and 
$B := \d ( L \cap \d S)$. 
Extend $X$ to $\X := X \sqcup \{\a \}$ by 
$\u \a := \{ \a \} \sqcup  A$ and 
$\d \a := \{ \a \} \sqcup  B$. 
Define a topology $\tau$ by an open subbase $\tau_X \cup \{\X -( \u E \cup \d F) \mid F, E \in X^\< \}$. 
Then $X$ is an embedded subspace which is not closed.
Thus it suffices to show the following claim, which induces that $X$ is not an H-closed pospace. 
\begin{claim}
$\X$ is a pospace. 
\end{claim}
Indeed, let $x \in X$ be any element. 
If $x \in A$, then $x >\a$. 
Since $\W_L(L \setminus \d S)$ does not exist,
there is an element $y$ of $L \setminus \d S$ such that $\a < y <x$. 
Then there is a finite subset $E$ of $A$ such that 
$y \in E$ and $E$ is a maximal chain. 
Since $(E- \{ y \} ) \subseteq I_y$, we have that   
$U:= \X -\d E$ is an open \nbd of $x$ and  $V := \X - \u E$ is an open \nbd of $\a$ such that 
$\d V \cap U = \emptyset $. 
If $x \in B$, then the symmetry of pospace implies that there are 
a finite subset $E$ of $X$ ,  
an open \nbd $U:= \X -\d E$ of $\a$ 
and an open \nbd $V := \X - \u E$ of $x$ such that 
$\d V \cap U = \emptyset $. 
Otherwise $x \| \a$. 
Since $L$ is a maximal chain, either 
$\u x \nsupseteq L \setminus \d S$ or   
$\d x \nsupseteq L \cap \d S$. 
First, consider the case $\u x \supseteq L \setminus \d S$.  
So $\d x \nsupseteq L \cap \d S$. 
Since $x \| \a$ and $\d x \nsupseteq L \cap \d S$, there is an element $\beta < \a $ such that $\beta \| x$. 
Let $E$ be a finite subset of $X$ such that $E \sqcup \{ x , \beta \}$ is a maximal chain. 
Since $\u x \supseteq L \setminus \d S$, 
we have $U := \X- (\u x  \cup \d x )$ is an open \nbd of $\a$ and 
$V := \X- (\u ( E \sqcup \{\beta \}) \cup \d ( E \sqcup \{\beta \}) ) \subseteq \u x \cup \d x$ is an open \nbd of $x$. 
Since $E \sqcup \{ x, \beta \}$ is a maximal antichain, 
Lemma \ref{lem:001} implies that 
$\u V \subseteq \u x \cup \d x$ and 
$\d V \subseteq \u x \cup \d x$.
Since $U \cap (\u x \cup \d x) = \emptyset$, we obtain 
$\u V \cap U = \emptyset$ and 
$\d V \cap U = \emptyset$. 
Second, consider the case the case $\u x \nsupseteq L \setminus \d S$.  
By symmetry, we may assume that 
$\d x \nsupseteq L \cap \d S$. 
Since $x \| \a$ and $\d x \nsupseteq L \setminus \d S$,
there is an element $\beta > \a $ such that 
$\beta \| x$.  
Since $x \| \a$ and $\d x \nsupseteq L \cap \d S$, 
there is an element $\gamma < \a $ such that 
$\gamma \| x$. 
Since $x \| \a$, there is a finite subset $E$ of $X$ such that 
$E \sqcup \{ x , \a \}$ is a maximal chain. 
Then 
$U := \X- (\u x \cup \d x )$ is an open \nbd of $\a$ and 
$V := \X- (\u ( E \sqcup \{\gamma \}) \cup \d ( E \sqcup \{\beta \}) ) \subseteq  (\u x \cup \d x) \setminus (\d \beta \cup \u \gamma)$ 
is an open \nbd of $x$. 
Since $E \sqcup \{ x,  \a \}$ is a maximal antichain and 
$V \subseteq \X- (\u ( E \sqcup \{ \a \}) \cup \d ( E \sqcup \{\a \}) )$, 
Lemma \ref{lem:001} implies that 
$\u V \subseteq \u x \cup \d x$ and 
$\d V \subseteq \u x \cup \d x$. 
Since $U \cap (\u x \cup \d x) = \emptyset$, we obtain 
$\u V \cap U = \emptyset$ and 
$\d V \cap U = \emptyset$. 

\end{proof}

\begin{lemma}\label{lem:004}
Let $X$ be an H-closed pospace without an infinite antichain. 
Then $X$ is directed complete and down-complete. 
\end{lemma}

\begin{proof}
Let $L$ be any infinite chain of $X$. 
Then $\min ( L^{\u}) = \min ( (\d L)^{\u})$. 
By Lemma \ref{lem:03}, we obtain 
$\min ( L^{\u}) \neq \emptyset$. 
Since $X$ has no infinite antichain, we have 
$\min ( L^{\u})$ is a nonempty finite subset. 
Put $\{ x_1, \dots , x_n \} = \min ( L^{\u})$. 
Since $X$ has no infinite antichain, there is a finite subset $K$ such that 
$K \supseteq \min ( L^{\u}) $ is a maximal antichain. 
If $\min ( L^{\u})$ is not a single set (i.e. $n >1$), 
then $x_i \| x_j$ for any distinct pair $i \neq j$. 
For any $i \neq j$, 
since $X$ is a pospace and $x_i \| x_j$ , 
there are open \nbds of $x_i , x_j$ respectively such that 
$\d U_i \cap U_j = \emptyset = U_i \cap \d U_j$. 
For any $i = 1, \dots , n$, there is an open \nbd  $U_i \subseteq 
(\u x_i \cup \d x_i) \setminus ( \u ( K - \{ x_i \} ) \cup \d ( K - \{ x_i \} ))
 $ of $x_i$ such that 
$U_i  \cap \d L =\emptyset$.   
Extend $X$ to $\X := X \sqcup \{ \a \}$ by 
$\u \a := \u_X \{ x_1, \dots , x_n \}  \sqcup \{ \a \}$ and 
$\d \a := \d_X L \sqcup \{ \a \}$
Define a topology $\tau$ on $\X$ by an open subbase 
$\tau_X \cup \{ \X - (\u F \cup \d E ) \mid F, E \in X^{\<} \}$, 
where $\tau_X$ is the topology of $X$. 
\begin{claim}
$\X$ is a pospace. 
\end{claim}

Indeed, let $x $ be any element of $X$. 
Let $G$ be a finite subset of $I_{\a}$ such that $G \sqcup \{ \a \}$ is a maximal antichain in $\X$.  
If $x = x_i$ for some $i$, 
then let $U := \X - (\u K \cup \u G \cup \d G ) \subseteq \d \a $ be an open \nbd of $\a$.  
Since $\d \a = \d_X L \sqcup \{ \a \}$,  
we have $\d U \cap U_i = \emptyset$. 
If $x \in \u \a - \{ x_1, \dots , x_n \}$, then let 
$U := \X - \u K $ be an open \nbd of $\a$ and 
$V:= \X - \d K \subseteq \u K$ an open \nbd of $\a$. 
Now $\d U \cap V = \emptyset$. 
If $x \in \d \a $, then there is an element $y \in L$ such that $x < y  < \a$. 
Now there is a finite subset $F$ of $X$ such that $y \in F$ and $F \cup \{ y \}$ is a maximal antichain. 
Then $F \subseteq I_y - \u y$. 
Hence $U:= \X - \d ( \{y \} \sqcup F)$ is an open  \nbd of $x$ and 
$V:= \X - \u ( \{y \} \sqcup F)$ is an open  \nbd of $x$ such that 
$U \cap \d V = \emptyset$. 
Otherwise $x \| \a$. 
Since $x \notin \d L$, there is an element $\beta < \a \in L$ such that 
$\beta \| x$. 
Then there is a finite set $F$ of $X$ such that $F \sqcup \{ x, \beta \}$ is a maximal antichain. 
Thus $U := \X- ( \d x \cup \u x )$ is an open \nbd of $\a$ and 
$V := \X- (\u ( F \sqcup \{\beta \}) \cup \d ( F \sqcup  \{ \beta \} ) ) \subseteq (\u x \cup \d x)$ 
is an open \nbd of $x$. 
Since $F \sqcup \{ x, \beta \}$ is a maximal antichain, 
Lemma \ref{lem:001} implies that 
$\u V \subseteq \u x \cup \d x$ and 
$\d V \subseteq \u x \cup \d x$.
Since $U \cup  (\u x \cup \d x) = \emptyset$, we obtain 
$\u V \cap U = \emptyset$ and 
$\d V \cap U = \emptyset$. 

Therefore $X$ is not H-closed, which contradicts. 
Thus $\min (  L^{\u})$ is a single set and so 
$\V L$ exists. 
\end{proof}

Notice that the symmetry of pospace implies that the dual statement of Lemma \ref{lem:004} holds 
(i.e. An H-closed pospace without an infinite antichain is down complete).


\begin{lemma}\label{lem:005}
Let $X$ be an H-closed pospace without an infinite antichain. 
For any nonempty chain $L \subseteq  X$, 
we have  $ \V L \in \cl \d L$.
\end{lemma}

\begin{proof}
Suppose that 
there is a chain $L \subseteq  X$ such that 
$ \V L \notin \cl \d L$. 
Put $a := \V L$. 
Let $\X := X \sqcup \{ \a \}$ be an extension of $X$ by 
$\u \a = \{ \a \} \sqcup \u a$ and 
$\d \a = \{ \a \} \sqcup \d L$. 
Define a topology $\tau$ on $\X$ by an open subbase 
$\tau_X \cup \{ \X - (\u F \cup \d E ) \mid F, E \in X^{\<} \}$, 
where $\tau_X$ is the topology of $X$. 
\begin{claim}
$\X$ is a pospace. 
\end{claim}

Indeed, let $x$ be an element of $X$. 
first we consider the case $x \in \u \a$. 
Since $X$ has no infinite antichain,  
there is a finite subset $F$ of $X$ such that 
$F \sqcup \{ a \} $ is a maximal antichain. 
Since $\u \a \cup \d \a \subseteq \u a \cup \d a$, we obtain 
$F \sqcup \{ \a \} $ is an antichain. 
If $x = a$, 
then $U = \X - ( \u a \cup \u F \cup \d F) \subseteq \d a$ is an open \nbd of $\a$. 
Since $ a \notin \cl \d L $, 
there is an open  \nbd $W$ of $a$ such that $W \cap \d L = \emptyset$. 
Then $V := W \setminus (\u F \cup \d F) \subseteq \d a \cup \u a$ is 
an open \nbd of $a$. 
Then $\d U \cap V \subseteq \d a - \d \a$. 
We may assume that $\d U \cap  V  \neq \emptyset$.  
Then $\d a - \d \a  \neq \emptyset$. 
Let $E$ be a maximal antichain of $\d a - (\d \a \sqcup \{ a \})$. 
Note that $E \cap (\u \a \cup \d \a) = \emptyset$. 
Since $X$ has no infinite antichain, we obtain $E$ is finite. 
Since $\a \notin \u E \cup \d E$, 
$U' := U \setminus ( \u E \cup \d E) $ is an open \nbd of $\a$. 
Since $\d a - \d \a \subseteq \d E \cup \u E$, we have 
$U' \subseteq \d \a$ and so 
$\d U' \cap V = \emptyset$. 
Otherwise $x >a$. 
$U = \X - \u ( F \cup \{ a \} ) $ is an open \nbd of $\a$ and  
$V := \X - \d (F \cup \{a \}) $ is an open \nbd of $x$ with $\d U \cap V = \emptyset$. 

Second, we consider the case $ x < \a$. 
Since $ a = \V L\notin L $, 
there are an element $y \in X$ and $ F \in X^\<$ such that 
$x < y < \a$ and $F \sqcup \{y \}$ is a maximal chain. 
Then $U := \X- \u ( F \sqcup \{y \})$ is an open \nbd of $x$ and 
$V := \X- \d ( F \sqcup \{y \})$ is an open \nbd of $\a$. 
Then $U \cap \d V = \emptyset$.

Finally, we consider the case $x \| \a$. 
If $x \| a$,  
then there is a finite subset $E$ of $X$ such that 
$E \sqcup \{ x , a \}$ is a maximal chain. 
Thus $U := \X- (\u x  \cup \d x) $ is an open \nbd of $\a$ and 
$V := \X- (\u ( E \sqcup \{a \}) \cup \d ( E \sqcup \{a \}) )$ is an open \nbd of $x$. 
Since $E \sqcup \{ x, a \}$ is a maximal antichain, 
Lemma \ref{lem:001} implies that 
$\u V \subseteq \u x \cup \d x$ and 
$\d V \subseteq \u x \cup \d x$.
Since $U \subseteq \d \a  \setminus (\u x \cup \d x)$, we obtain 
$\u V \cap U = \emptyset$ and 
$\d V \cap U = \emptyset$. 
Otherwise $x < a$. 
Since $x \| \a$ and 
$a = \V L \notin L$, 
there is an element $\beta \in L$ such that $\beta \| x$. 
Let $E$ be a finite subset of $X$ such that 
$E \sqcup \{ x , \beta \}$ is a maximal antichain. 
Then $U := \X - ( \d x \cup \u x ) $ is an open \nbd of $\a$ and 
$V := \X - ( \d (E \sqcup \{ \beta \} ) \cup \u (E \sqcup \{ \beta \} )) $ 
is an open \nbd of $x$. 
Since $E \sqcup \{ x, \beta \}$ is a maximal antichain, 
Lemma \ref{lem:001} implies that 
$\u V \subseteq \u x \cup \d x$ and 
$\d V \subseteq \u x \cup \d x$.
Since $U \subseteq \d \a  \setminus (\u x \cup \d x)$, we obtain 
$\u V \cap U = \emptyset$ and 
$\d V \cap U = \emptyset$. 
This completes the proof of this claim. 

By the definition of $\tau$, $X$ is a proper subspace and $\cl X = \X$. 
Therefore $X$ is not an H-closed pospace. 
\end{proof}

Notice that the symmetry of pospace implies that the dual statement of Lemma \ref{lem:004} holds 
(i.e. An H-closed pospace without an infinite antichain is down complete).

Now we show the another direction. 

\begin{lemma}\label{lem:05}
Let $X$ be a pospace without an infinite antichain. 
Suppose that  $X$ is directed complete and down-complete such that 
$ \V L \in \cl \d L$  
and $ \W L \in \cl \d L$ for any nonempty chain $L \subseteq  X$. 
Then $X$ is an H-closed pospace. 
\end{lemma}

\begin{proof}
Suppose there is a non-H-closed pospace $X$ without an infinite antichain such that 
$X$ is directed complete  and down-complete, and  
$ \V L \in \cl \d L$ and $ \W L \in \cl \d L$ for any chain $L \subseteq  X$. 
Then there is an embedding from $X$ to a pospace $\X$ such that 
$X$ is a dense proper subspace. 
Fix any element $x \in \X -X$. 
Since $X$ is directed complete and down-complete , 
$X = \d_X \max_X X = \u_X \min_X X$. 
Since $X$ has no infinite antichain, 
we have that $\max_X X$ and $\min_X X$ are finite subsets. 
If all elements of $X$ are incomparable to $x$, 
then $x \notin \d \max_X X$. 
Since $\d \max_X X$ is $\X$-closed, we have that 
$x \notin \cl X$, which contradicts to the density of $X$. 
Thus there is a comparable element $\a$ of $X$ to $x$. 
By the symmetry of pospace, 
we may assume that $\a <x$. 
If $x \notin \d \max_X X$, 
then the fact $\d \max_X X \supseteq X$ implies that 
$x \notin \cl X$, which is impossible. 
Thus there is an element $\omega \in X$ such that $x < \omega$. 
Since $X$ is embedded into $\X$, 
we have that 
$A' := \d x \cap X$ and 
$B' := \u x \cap X$ are $X$-closed. 
Let 
$A := \{ \V_X L \mid L \neq \emptyset \subseteq A' \text{ is a chain} \}$ 
and 
$B := \{ \W_X L \mid L \neq \emptyset \subseteq B' \text{ is a chain} \}$. 
Since $X$ is directed complete, 
we obtain 
$\max_X A $ and $\min_X B$ are nonempty such that   
$\d_X \max_X A \supseteq A$ and 
$\u_X \max_X B \supseteq B$. 

\begin{claim}
$x \notin \d \max_X A$. 
\end{claim}
Indeed, 
suppose that 
there is an element $y \in \max_X A$ such that $x < y$. 
By the definition of $A'$, we have $y \notin A'$. 
Hence there is a chain $L \subseteq A'$ such that 
$y = \V_X L \notin L$. 
By the assumption, 
$y \in \cl _X L \subseteq \cl _X A' =A'$, which is impossible.

By the symmetry of pospace, 
$x \notin \u \min_X B$. 
Let $I := \{ y \in X \mid x \| y \}$. 
Since $X$ has no infinite antichain, 
there is a finite subset $F$ of $I$ such that 
$I \subseteq \u F \cup \d F$. 
Since $X = A' \sqcup B' \sqcup I$, 
$\d_X \max_X A \supseteq A'$ and 
$\u_X \max_X B \supseteq B'$, 
we have $X \subseteq \d \max_X A \cup \u \min_X B \cup \u F \cup \d F$. 
Since $x \notin  \d \max_X A \cup \u \min_X B \cup \u F \cup \d F$, 
we obtain that 
$x \notin \cl X$, which is impossible. 
\end{proof}

Theorem \ref{thm:001} is induced by Lemma \ref{lem:004}, \ref{lem:005}, their dual statements, and  Lemma \ref{lem:05}.


\begin{thebibliography}{00}
\bibitem[GPR]{GPR}
O. Gutik-D. Pagon-D. Repov\v s
\textit{On Chains in H-Closed Topological Pospaces}
Order (2010) 27: 69-81
\bibitem[GR]{GR}
O. Gutik-D. Repov\v s
\textit{On linearly ordered H-closed topological semilattices} 
Semigroup Forum (2008) 77: 474-481
\bibitem[W]{W}
L. E. Ward, Jr.
\textit{Partially ordered topological spaces}
Proc. Amer. Math. Soc. (1954) 5:1, 144-161.
\end{thebibliography}
\end{document}